\documentclass[12pt]{article}

\usepackage{amsmath}
\usepackage{setspace}
\usepackage{array}
\usepackage{times}
\usepackage{bm}
\usepackage{graphicx}
\usepackage{latexsym}
\usepackage{color}

\usepackage[dvips, letterpaper, nohead, top=1.3in, bottom=1.3in, left=1.3in, right=1.3in]{geometry}

\newtheorem{theorem}{Theorem}
\newtheorem{lemma}{Lemma}

\newtheorem{remark}{Remark}
\newtheorem{example}{Example}

\newtheorem{procedure}{Procedure}
\def\qed{\rule{2mm}{2mm}}

\normalsize

\begin{document}

\small\normalsize

\title{On Stepwise Control of Directional Errors under Independence and Some Dependence}
\author
{Wenge Guo \\
Department of Mathematical Sciences\\
New Jersey Institute of Technology \\
Newark, NJ 07102-1982 \\
\and
Joseph P. Romano \\
Departments of Statistics and Economics\\
 Stanford University\\
Stanford, CA 94305-4065\\
}

\date{}

\maketitle

\begin{abstract}
In this paper, the problem of error control of stepwise multiple testing procedures is considered.
For two-sided hypotheses, control of both type 1 and type 3 (or directional) errors is required,
and thus mixed directional familywise error rate control  and mixed directional false discovery rate control  are each considered  by incorporating both types of errors in the error rate.
Mixed directional familywise error rate  control of stepwise methods  in multiple testing has proven to be a challenging problem, as  demonstrated in Shaffer (1980).  By an appropriate formulation of the problem, some new stepwise procedures are developed that  control type 1 and directional errors under independence and various dependencies.
\end{abstract}

\normalsize

\section{Introduction}

The main problem  considered in this paper  is the construction of procedures for the  simultaneous testing of $n$ parameters $\theta_i$.
For convenience, the null hypotheses $\theta_i = 0$ are of interest.
 Of course, we would like to reject any
null hypothesis if the data suitably dictates, but  we also  wish to
make directional inferences about the signs of $\theta_i$.
First,  consider the problem of simultaneously testing  $n$ null hypotheses against two-sided alternatives:
\begin{equation}\label{equation:no1}
 \check{H}_i: \theta_i=0~~{vs.~~} \check{H}'_i: \theta_i \neq 0,~~ i = 1,
\ldots, n~.
\end{equation}
Suppose, for $i = 1, \ldots, n$,  a test statistic
 $T_i$, is available for testing $\check{H}_i$.
  If
$\check{H}_i$ is rejected, the decision regarding $\theta_i > 0$ (or
$\theta_i < 0$) is made by checking if $T_i > 0$ (or $T_i <0$). In
making such rejection and directional decisions, three types of
errors might occur. The first one is the usual type 1 error, which
occurs when $\theta_i = 0$, but we falsely reject $\check{H}_i$ and
declare $\theta_i \neq 0$. The second one is the type 2 error, which
occurs when $\theta_i \neq 0$, but we fail to reject $\check{H}_i$.
The last one is called type 3 or directional error, which occurs
when $\theta_i > 0$ (or $\theta_i < 0)$, but we falsely declare
$\theta_i < 0$ (or $\theta_i > 0)$. We wish to control both type 1
and type 3 errors at pre-specified levels and, subject to their control, find testing methods with small probability of type 2
errors.

Given any procedure which makes rejections as well as directional claims about any rejected hypotheses,  let  $\check{V}$ and $\check{S}$ denote the  numbers of
type 1 errors and type 3 errors, respectively,  among $\check{R}$ rejected
hypotheses.  Let  $\check{U} = \check{V} + \check{S}$ denoting the
total number of type 1 and type 3 errors.  Then, the usual familywise error rate (FWER) and false discovery rate (FDR) are
defined respectively by FWER = Pr$(\check{V} \ge 1)$ and FDR = $E
\left ( \check{V}/\max(\check{R}, 1) \right )$, and the mixed
directional FWER and FDR are defined respectively by mdFWER =
Pr$(\check{U} \ge 1)$ and mdFDR = $E \left (
\check{U}/\max(\check{R}, 1) \right )$.

The main objective of this paper is to develop  stepwise
procedures (described shortly)  for controlling the mdFWER and mdFDR when simultaneously
testing the $n$ two-sided hypotheses $\check{H}_1, \ldots,
\check{H}_n$.
In multiple testing, the problem of simultaneously
testing $n$ two-sided hypotheses along with directional decisions
subject to the control of the mdFWER is technically very
 challenging. Until now, only a few results have been  obtained
under the strong assumption of independence of the test statistics
along with some additional conditions on the marginal distribution
of the test statistics.

Shaffer (1980) proved that if the test statistics $T_i, i = 1,
\ldots, n$ are mutually independent and if the distributions of the
$T_i$'s satisfy some additional conditions, the mdFWER of a
directional Holm procedure is strongly controlled at level $\alpha$.
She also constructed a counterexample where the aforementioned
procedure loses the control of the mdFWER even under independence
when the test statistics are Cauchy distributed. Holm (1979b, 1981)
extended Shaffer's  (1980)'s result to normal distributional settings
where the $T_i$'s are conditionally independent. Finner (1994) and
Liu (1997) independently used Shaffer's  (1980) method of proof to
show the mdFWER control of directional Hochberg procedure by
making the same distributional assumptions  as Shaffer (1980). By
generalizing Shaffer's method of proof, Finner (1999) extended
Shaffer's result on the  Holm procedure to a large class of stepwise or
closed multiple testing procedures under the same assumptions as in
Shaffer (1980). He also gave a new but very simple and elegant proof
for the aforementioned result under the assumption of TP$_3$
densities. For further discussions on the mdFWER control of closed testing methods, see
Westfall, Bretz and Tobias (2013).

Another method to tackle the problem of directional errors has been
considered in Bauer, Hackle, Hommel and Sonnemann (1986), in which the
problem of  testing $n$ two-sided hypotheses testing with additional
directional decisions is reformulated as the problem of  testing $n$ pairs
of one-sided hypotheses given by $$H_{i1}: \theta_i \le 0~~{\rm  vs.}~~
H_{i1}^{'}: \theta_i > 0~, \nonumber$$ and
 $$\tilde{H}_{i2}: \theta_i \ge 0~~ {\rm vs.~~}
\tilde{H}_{i2}^{'}: \theta_i < 0~\nonumber$$ for $i=1, \ldots, n$. They proved that
without additional distributional assumptions, only a slight
improvement of the conventional Holm procedure is possible for
testing these $2n$ hypotheses. They also showed by a counterexample
that in general distributional settings, a further improvement of
their procedure is impossible. Compared with Shaffer's  (1980)
directional Holm procedure for testing $n$ two-sided hypotheses,
their procedure is very conservative, although it controls
directional errors under more general distributional settings of
arbitrary dependence.

Finally, they also reformulated the
aforementioned problem as the  problem of  testing $n$ pairs of
one-sided hypotheses given by $$H_{i1}: \theta_i \le 0~~{\rm vs.~~}
H_{i1}^{'}: \theta_i > 0~,\nonumber$$ and
 $$H_{i2}: \theta_i > 0 ~~{\rm  vs.~~}  H_{i2}^{'}:
\theta_i < 0~,\nonumber$$
for $i=1, \ldots, n$, among which there is exactly
one true null hypothesis within each pair of one-sided hypotheses.
They proved that the modified Bonferroni procedure with the critical
constant $\alpha/n$ (as opposed to $\alpha/ 2n$) strongly controls the FWER when testing these
$2n$ one-sided hypotheses.  This result is of course trivial because in this formulation
there are exactly $n$ true null hypotheses.  At the same time, given that there are always $n$ true null hypotheses, it is perhaps surprising that one can, as we do, develop stepdown methods
that improve upon this single step method.  (Indeed, at any step when applying a stepdown method,
there are always $n$ true null hypotheses, and this number does not reduce.)

 In the above two formulations of one-sided hypotheses, there are some inherent
disadvantages when developing stepwise methods for controlling
the FWER. In the first formulation, there may be a different number of
true null hypotheses between $\theta_i = 0$ and $\theta_i \neq 0$, which makes it challenging to develop powerful stepwise methods in this formulation, as shown  in Bauer et al. (1986).
In the second formulation, one possible type 1 error will not be counted even though
$T_i$ is very small when $\theta_i = 0$, which makes it unable to completely control type 1
and type 3 errors in the original formulation of two-sided hypotheses even though the FWER is controlled in this formulation.  Further discussion of this point will be presented later.
On the other hand, the problem of the mdFDR control seems to be technically less challenging and methods for controlling the mdFDR are available (see Benjamini and Yekutieli, 2005; Guo, Sarkar and Peddada, 2010).

In the next section, some basic notation is given, as well as our approach to the problem.
Theorems 1--4 deal with control of the familywise error rate with directional decisions,
first under independence, and then under block dependence and positive dependence.
Theorems 5--8 analogously provide results for the false discovery rate.

Although many procedures are introduced in this paper, their proven control of the FWER or FDR
are established under different assumptions of dependence, including independence,  between-block dependence, within-block dependence, and positive dependence. It would be impossible to advocate
a single procedure in applications without any knowledge of dependence. It would be more appropriate to suggest different procedures based on different dependence information. Only under the same assumption of independence, four different procedures, Procedures 1-3 and Procedure $1'$, are developed for controlling the FWER. Among them, we recommend
the use of Procedure 3 in practice because this procedure is generally the most powerful
while controlling the FWER in the sense that its critical values are generally larger. The main
reason for introducing Procedures 1, 1', and 2 were as building blocks to the stepdown method
of Procedure 3. Procedure 4 is developed in order to control the FWER under block dependence
(Theorem 3), while Procedure 5 applies to positive dependence (Theorem 4). Procedure
6 applies to control of the FDR under independence (Theorem 5). Procedure 7 applies to
control of the FDR under between-block dependence (Theorem 6), while Procedure 8 applies
to control of the FDR under within-block dependence (Theorem 7). Procedure 9 applies to
FDR control under positive dependence (Theorem 8).

\section{Preliminaries}

In this section,  some necessary notation and basic concepts are introduced.

\subsection{Notation}

Suppose $T_i$ has  cumulative distribution function $F_{i, \theta_i}(\cdot)$ (with density denoted
$f_{i, \theta_i}(\cdot)$ when it is well-defined),
both of which depend on a single parameter $\theta_i$.
It is assumed that the null distribution of $T_i$, i.e. $F_{i, 0} ( \cdot )$ is continuous.
We also assume that $F_{i, \theta_i}(t)$ is
non-increasing in $\theta_i$ for any given $t$ and $F_{i, 0}(t)$
is symmetric about zero, i.e., $F_{i, 0}(-t) = 1 - F_{i, 0}(t)$ for any $t$.
(In fact, the symmetry assumption is not really necessary; indeed, one may take the probability integral transformation $F_{i,0} ( T_i )$  to get
a new test statistic  that is uniform and then shift it by 1/2 to get a ``symmetric" null test statistic.)
Let $t_i$ be the observed  value of $T_i$. Then, the (two-sided)  $p$-value for
testing $\check{H}_i$ is
$$\check{P}_i = 2 \min(F_{i, 0}(t_i), 1 -
F_{i, 0}(t_i))~.$$
Let $\check{P}_{(1)} \le \ldots \le \check{P}_{(n)}$ be
the ordered $p$-values and $\check{H}_{(1)}, \ldots,
\check{H}_{(n)}$ the associated null hypotheses. Then, given a
non-decreasing set of critical constants $0 < \alpha_1 \le \cdots
\le \alpha_n < 1$, a stepdown multiple testing procedure rejects the
set of null hypotheses $\{\check{H}_{(i)}, i \le i_{SD}^*\}$ and
accepts the rest, where $i_{SD}^* = \max \{i: \check{P}_{(j)} \le
\alpha_j~ \forall~j \le i \}$ if the maximum exists, and  otherwise it
accepts all the null hypotheses. A stepup procedure, on the other
hand, rejects the set $\{\check{H}_{(i)}, i \le i_{SU}^*\}$ and
accepts the rest, where $i_{SU}^* = \max \{i: \check{P}_{(i)} \le
\alpha_i \}$ if the maximum exists, otherwise it accepts all the
null hypotheses. Furthermore, if stepwise procedures (stepdown or
stepup) are applied along with additional directional decisions,
such procedures are often termed as directional stepwise procedures
(Shaffer, 2002). (A stepwise procedure with constant $\alpha_i$ is referred to
as a single-step procedure.)
The constants in a stepwise procedure are determined subject to the
control of a suitable error rate at a pre-specified level $\alpha$.

\subsection{Formulation}

In order to further explore the problem of controlling type 1 and type 3 errors under
independence, and also under some
dependence, we  first reformulate this problem as an equivalent  one of
simultaneously testing multiple one-sided hypotheses subject to the
control of the FWER (or FDR). Specifically, $\check{H}_i, i= 1,
\ldots, n$ against two-sided alternative $\check{H}'_i$ is reformulated as three
null hypotheses $H_{ij}, j = 1, 2, 3$ against one-sided alternatives $H'_{ij}$,
$$H_{i1}:
\theta_i \le 0~~{\rm  vs.}~~ H_{i1}^{'}: \theta_i > 0~, \nonumber $$
$$H_{i2}: \theta_i >
0~~{\rm vs.}~~H_{i2}^{'}: \theta_i < 0~,\nonumber $$
 and $$H_{i3}: \theta_i = 0~~{\rm  vs.}~~
H_{i3}^{'}: \theta_i < 0~.\nonumber$$
 As we know, for the original problem of testing the
 two-sided hypotheses $\check{H}_i, i= 1, \ldots, n$ along with
directional decisions, there are two possibilities of type 1 errors
and two possibilities of type 3 errors.   Indeed, when $\theta_i =
0$, the corresponding test statistic $T_i$  can be too  large or too
small;  or,  when $\theta_i
> 0~ (\text{or } < 0)$, $T_i$ is too small (or large). In the new
formulation, those two possible directional errors in the original
problem are transformed as type 1 errors for testing $H_{i1}$ and
$H_{i2}$, respectively, and the two possible type 1 errors when
testing $\check{H}_i$ are transformed as type 1 errors for testing
$H_{i1}$ and $H_{i3}$, respectively.
 It should be noted that
the additional directional decisions in all these formulations of one-sided alternatives is
unnecessary as any rejection already corresponds to a directional decision.
Note that,
when $T_i$ is used for testing  $H_{i1}$,
$-T_{i}$ is used for testing both $H_{i2}$ and $H_{i3}$.

Let $\mathcal{F} = \{H_{ij}: i = 1, \ldots, n, j = 1, 2, 3\}$ denote
the whole family of the $3n$ hypotheses $H_{ij}$'s to be
tested. We split $\mathcal{F}$ as two subfamilies $\mathcal{F}_1$
and $\mathcal{F}_2$, where
$$\mathcal{F}_1 = \{H_{ij}: i = 1, \ldots,
n, j = 1, 2\} \nonumber $$ and $$\mathcal{F}_2 = \{H_{i3}: i = 1, \ldots, n \}~.\nonumber$$
In this paper, we use a separate approach for testing multiple
families of hypotheses. In this approach, two given
multiple testing methods are used for testing $\mathcal{F}_i, i = 1,
2$, respectively. If $S_i, i = 1, 2$ denote the respective
rejection sets for testing $\mathcal{F}_i$, then the rejection set
for testing $\mathcal{F} = \mathcal{F}_1 \bigcup \mathcal{F}_2$ is $S_1 \bigcup S_2$.
The advantage of splitting $\mathcal{F}$ derives from the fact
that $\mathcal{F}_1$ consists of $2n$ hypotheses, of which exactly $n$
of them correspond to true null hypotheses.

For the aforementioned approach, let $V_i$ denote the number of type
1 errors among $R_i$ rejected hypotheses when testing
$\mathcal{F}_i$ for $i = 1, 2$, and let $V$ denote the number of
type 1 errors among $R$ rejected hypotheses when testing
$\mathcal{F}$. Thus, $R = R_1 + R_2$ and $V = V_1 + V_2$. Then, the
FWER and FDR of the multiple testing method for testing
$\mathcal{F}$ are defined respectively by
$${\rm FWER}_{\mathcal{F}} =
Pr \{ V \ge 1 \} \nonumber$$
 and
 $${\rm FDR}_{\mathcal{F}} = E \left (V/\max(R, 1)
\right )~.\nonumber$$
 Similarly, the FWER and FDR for  testing $\mathcal{F}_i$
are defined respectively by FWER$_{\mathcal{F}_i}$ = Pr$ \{ V_i \ge 1 \}$
and FDR$_{\mathcal{F}_i}$ = $E \left (V_i/\max(R_i, 1) \right ), i =
1, 2$. Note that $V \ge 1$ implies $V_1 \ge 1$ or $V_2 \ge 1$,
so that  $\text{FWER}_{\mathcal{F}} \le \text{FWER}_{\mathcal{F}_1} +
\text{FWER}_{\mathcal{F}_2}.$ Similarly, using the simple  inequality
 $$\frac{V}{\max(R, 1)} = \frac{V_1}{\max(R, 1)} +
\frac{V_2}{\max(R, 1)} \le \frac{V_1}{\max(R_1, 1)} +
\frac{V_2}{\max(R_2, 1)}~,$$ we have $\text{FDR}_{\mathcal{F}} \le
\text{FDR}_{\mathcal{F}_1} + \text{FDR}_{\mathcal{F}_2}.$ We will
develop in this paper respective stepwise methods for controlling
the $\text{FWER}_{\mathcal{F}}$ and $\text{FDR}_{\mathcal{F}}$ when
testing $\mathcal{F}$ based on the aforementioned separate approach
and the above two inequalities. We note that in the existing
literature, a number of powerful stepwise methods have been introduced under
various dependencies for testing $\mathcal{F}_2$, for which unlike $\mathcal{F}_1$, there is
no specific dependency relationship among the test statistics
corresponding to those hypotheses in $\mathcal{F}_2$.
For example,   control of the
 $\text{FWER}_{\mathcal{F}_2}$
can be done by   the Holm (1979a)  and Hochberg (1988)
while   the Benjamini and Hochberg (1995) procedure (BH) can be used to control
 the $\text{FDR}_{\mathcal{F}_2}$ (Therefore, through much of the  paper, we will focus
on developing  stepwise methods for controlling the
$\text{FWER}_{\mathcal{F}_1}$ and $\text{FDR}_{\mathcal{F}_1}$ under
independence and certain dependencies, unless noted otherwise.)

Before we embark upon control of  any error rate for $\mathcal{F}_1$ as a building block
for control over the larger family $\mathcal{F}$, we would like to argue that
this seemingly more restrictive control over the smaller family $\mathcal{F}_1$ is already a plausible
approach to the problem of control of  directional errors.
For this, we draw upon the wisdom and philosophy of  one of the fathers in the field of
multiple testing, John Tukey.    In the context of single testing, Tukey argued
that a point null hypothesis is never true, and therefore control of type 1 errors is
the wrong formulation.  Tukey cared more about whether or not one could tell the ``effect size"
or the ``sign" of a parameter.  To quote Tukey (1991), ``Statisticians  classically asked the wrong question -- and were willing to answer with a lie, one that was often a downright lie.......All we know
about the world teaches us that the effects of $A$ and $B$ are always different -- in some decimal place -- for any $A$ and $B$.  Thus asking `Are the effects different' is foolish.   What we should
be answering first is `Can we tell the direction in which the effects of $A$ differ from the effects
of $B$?'. "  Thus, for Tukey, emphasis must be completely upon control of directional or type 3 errors.
So, as also argued in Jones and Tukey (2000) in the context of a test of a single parameter
$\theta$ (which is motivated there as a difference in means), one can and should apply
a classical two-sided $t$-test  so that the probability of observing an outcome in either
the right or left tail is not $\alpha /2$, but $\alpha$.  That is,   if one wishes  to make directional
inferences or claims about a parameter (which is always desirable)  then the problem of testing  the null  hypothesis $\theta = 0$
 at level $\alpha$ should be replaced by the problem of testing the two hypotheses:  testing
$ \theta < 0$ against $\theta > 0$ as well as testing $\theta > 0$ against $\theta < 0$.
 Since $\theta = 0$ never holds, one can always use the $1- \alpha$ quantile in the right tail
 rather than the $1 - \alpha/2$ quantile, and similarly the $\alpha$ quantile in the left tail.
 In our context, if we acknowledge that $\theta_i$ is never 0 from the start, then
 we never need to include $\mathcal{F}_2$ in the family of hypothesis tested, and the problem
 of control of directional errors is equivalent to control of the error rate over $\mathcal{F}_1$.
 Moreover, if one takes Tukey's stance  to heart, then the inequality in the definition
 of $H_{i1}$ can be a strict inequality.   However, we retain the inequality
 because the methods we develop apply to $H_{i1}$ as defined, and hence to the more restricted definition.  Thus,  control over $\mathcal{F}_1$ is emphasized throughout, as both a building
 block toward control over $\mathcal{F}$ but also as a formulation worth studying in its own right.
 A nice review of Tukey's contributions to multiple testing can be found in Benjamini and Braun (2002).

\subsection{Assumptions}

It should be noted that   $H_{i1} \bigcap
H_{i2}$ is empty and $H_{i1} \bigcup H_{i2}$ is the whole parameter space.
 Thus, there are exactly $n$ true and $n$ false null
hypotheses in $\mathcal{F}_1 = \{H_{ij}: i = 1, \ldots,
n, j = 1, 2\}$, which form $n$ pairs of true and false nulls $(H_{i1}, H_{i2})$.
For notational convenience, we
respectively use $H_1, \ldots, H_n$  and $H_{n+1}, \ldots, H_{2n}$
denoting the $n$ true and $n$ false nulls with $(H_i, H_{n+i})$ denoting
$(H_{i1}, H_{i2})$ and $(P_i,
P_{n+i})$ denoting the pair of the corresponding (one-sided) $p$-values. With
the test statistic $T_i$ and the calculated value $t_i$, the
$p$-value $P_i$ corresponding to $H_i$ is equal to $F_{i, 0}(t_i)$ or $1
- F_{i, 0}(t_i)$ depending on $H_i: \theta_i \le 0$ or $H_i: \theta_i >
0$, and $P_{n+i} = 1 - P_i$ for $i = 1, \ldots, n$. In addition, let
 $I_0 = \{1, \ldots, n\}$ and $I_1 = \{n+1, \ldots, 2n\}$
denote the index sets of true and false nulls among the $2n$
hypotheses, $H_1, \ldots, H_{2n}$, respectively.

Regarding the marginal distribution of the true null $p$-values,
 the following  assumptions are invoked  throughout much of the  paper:

\vskip 5pt

\textbf{A.1 } For any $p$-value $P_i, i \in I_0$ and given parameter
$\theta_i$,
\begin{eqnarray}\label{eq:1}
\text{Pr}_{\theta_i} \left \{P_i \le p \right \} \le p \text{~~ for~any~ } 0 \le p \le
1.
\end{eqnarray}
For $\theta_i = 0$, (\ref{eq:1}) is an equality;  that is, $P_i \sim
U(0, 1)$ for $i \in I_0$ when $\theta_i = 0$.

\textbf{A.2 } For any $p$-value $P_i, i \in I_0$ and given parameter
$\theta_i$,
\begin{eqnarray}\label{eq:2}
\text{Pr}_{\theta_i} \left \{ P_i \le p \big | P_i \le p' \right \}
\le \text{Pr}_{\theta_i = 0} \left \{ P_i \le p \big | P_i \le p'
\right \} ~.
\end{eqnarray}
for any $0 \le p \le p' \le 1.$

\textbf{A.3 } The test statistics $T_i, i = 1, \ldots, n$ are
mutually independent.

\vskip 10pt

While the assumption of independence is quite restrictive, to the best of
our knowledge, all the previous results on the mdFWER control of the existing stepwise
procedures along with directional decisions are established under this assumption.
However, not all of our  results require both A2 and A3.

Of course, under assumption A.1, the right hand side of (\ref{eq:2}) is just $ p / p'$.
Assumption A.2 is easily satisfied by the usual test statistics.
Actually,  the following result holds.

\begin{lemma}\label{lem:1}
If  the family of densities $f_{i, \theta}(\cdot)$  of $T_i$ satisfies the assumption of
monotone likelihood ratio (MLR), i.e., for any given $\theta_1 >
\theta_0$ and $x_1
> x_0$,
$ \frac{f_{i, \theta_1}(x_1)}{f_{i, \theta_0}(x_1)} \ge
\frac{f_{i, \theta_1}(x_0)}{f_{i, \theta_0}(x_0)}$, then Assumption A.2 holds.
\end{lemma}
For the proof of Lemma 1, see the  Appendix.  Of course, the assumption holds if the distribution of $T_i$ is a normal shift model, which often asymptotically approximates the underlying situation.

By Lemma \ref{lem:1},  the MLR assumption implies Assumption A.2.
However, these two assumptions are not equivalent. Assumption A.2 is
slightly weaker than the MLR assumption.  It is equivalent to the
following condition: for any given $\theta_1$ and $x_1
> x_0$, $\frac{F_{i, \theta_1}(x_1)}{F_{i, 0}(x_1)} \ge
\frac{F_{i, \theta_1}(x_0)}{F_{i, 0}(x_0)}$ when $\theta_1 > 0$ and
$\frac{1 - F_{i, \theta_1}(x_1)}{1-F_{i, 0}(x_1)} \le \frac{1-
F_{i, \theta_1}(x_0)}{1-F_{i, 0}(x_0)}$ when $\theta_1 < 0$. It should be
pointed out that Assumption A.2 is different from the conventional
TP$_2$-property of $\partial[1 - F_{i, \theta_i}(x)]/{\partial
\theta_i}$, which is almost  always assumed in  the existing
literature on control of directional errors (Shaffer, 1980; Finner,
1999). The only exception is Sarkar, Sen and Finner (2004). In that
paper, it is assumed that  $f_{i, \theta_i}(\cdot)$   satisfies  the
aforementioned MLR condition.

To characterize the joint distribution among the test statistics
$T_i, i = 1, \ldots, n$, several dependence assumptions have been
made in this paper: independence, within-block dependence,
between-block dependence, and positive dependence. The positive
dependence condition, which will be of the type characterized by the
following:
\begin {eqnarray}\label{eq:3}
E \left \{ \phi(T_1, \ldots, T_n) ~|~ T_i \ge u \right \} \uparrow u
\in (0,1),
\end {eqnarray}
for each $T_i$ and any (coordinatewise) non-decreasing function
$\phi$. This type of positive dependence is commonly encountered and
used in multiple testing; see, for instance, Sarkar (2008) for
references. Other dependence conditions such as independence,
within-block and between-block dependence, will be characterized in
Sections 3 and 4, respectively.

\section{Controlling the mdFWER under independence}

In this section,  several stepwise procedures for
controlling the $\text{FWER}_{\mathcal{F}_1}$ are presented  under the
assumption of independence. \vskip 5pt

\subsection{Two-stage procedure}

For simplicity, we first consider a two-stage version of the usual
Holm procedure for testing $\mathcal{F}_1$ as
follows.
\begin{procedure}
(Two-stage procedure) \rm
\begin{enumerate}
  \item Reject all null hypotheses $H_i$ with the $p$-values
less than or equal to $\alpha/n$. Let $r$ be the total number of
rejections at this stage. If $r = n$, we stop testing; otherwise,
  \item For the remaining hypotheses, reject those with the
$p$-values less than or equal to $\alpha/(n - r)$.
\end{enumerate}
\end{procedure}

In the above Procedure 1, the Bonferroni procedure is used in the first stage
for testing the $2n$ hypotheses. Generally, the
Bonferroni would actually use the critical constant $\alpha/2n$ when
testing $\mathcal{F}_1$. However, in this formulation we know there
are exactly $n$ true null hypotheses in $\mathcal{F}_1$ and we can
apply an obviously modified Bonferroni procedure with critical
constant $\alpha/n$. Our method then improves upon this with a
second stage improvement in the spirit of a stepdown method.
Procedure 1 can also be regarded as an adaptive
Bonferroni procedure with the critical constant $c = \alpha/\max(n -
R_1(\alpha/n), 1)$, where $R_1(\alpha/n) = \sum_{i=1}^{2n} I(P_i \le
\alpha/n)$ (Finner and Gontscharuk, 2009; Guo, 2009).

For any given parameter
vector $\theta = (\theta_1, \ldots, \theta_n)$, we have
\begin{eqnarray}\label{eq:6}
\text{FWER}_{\mathcal{F}_1}(\theta) & \le &
\frac{\alpha}{1 - \alpha/n},
\end{eqnarray}
whose proof is given in the Appendix.

\begin{theorem}\label{thm:1}
Consider Procedure 1 defined as above. Under assumptions
A.1 - A.3, the following conclusions hold.
\begin{enumerate}
  \item[\emph{(i)}] The procedure strongly controls the
  $\text{FWER}_{\mathcal{F}_1}$ at level $\frac{\alpha}{1 -
  \alpha/n}$.
  \item[\emph{(ii)}] $\limsup_{n \rightarrow \infty}
  \text{FWER}_{\mathcal{F}_1}\le \alpha.$ That is, the procedure asymptotically controls the
  $\text{FWER}_{\mathcal{F}_1}$ at level $\alpha$. Moreover,
  if the critical constants of the two-stage directional procedure are rescaled by using
$\frac{\alpha}{1 + \alpha/n}$ to replace $\alpha$, then the
resulting procedure, which is labeled as Procedure $1'$,  strongly controls the
$\text{FWER}_{\mathcal{F}_1}$ at level $\alpha$ even in  finite
samples.
\end{enumerate}
\end{theorem}

\begin{remark}\rm

It should be noted that Procedure $1'$ in
Theorem \ref{thm:1} is not consistently more powerful than Bauer et al.
(1986)'s modified Bonferroni procedure with the critical constant
$\alpha/n$, since its critical constant at stage 1 is slightly
smaller than $\alpha/n$. However, by carefully checking the proof of
(\ref{eq:6}) (see the Appendix), we can see that for Procedure 1, we actually only need to rescale its critical constant at
stage 2 in order to maintain the control of the FWER at level
$\alpha$. The newly modified procedure is described in details as
follows.
\begin{procedure}
(Modified two-stage procedure) \rm
\begin{enumerate}
  \item Reject all null hypotheses $H_i$ with the $p$-values
less than or equal to $\alpha/n$. Let $r$ be the total number of
rejections at this stage. If $r = n$, we stop testing; otherwise,
  \item For the remaining hypotheses, reject those with the
$p$-values less than or equal to $\beta/(n - r)$, where $\beta =
\frac{\alpha}{1 + \alpha/n}$.
\end{enumerate}
\end{procedure}

It is easy to see that the above Procedure 2 is consistently more
powerful than Bauer et al.'s modified Bonferroni procedure, because
for this procedure, even if only one hypothesis is rejected at
stage 1, its critical constant $\frac{\alpha}{(n-1)(1 + \alpha/n)}$
at stage 2 is also larger than $\alpha/n$, the critical constant of
Bauer et al.'s procedure.

\end{remark}

\begin{remark}\rm
Goeman and Solari (2010) recently provided a very general approach for
developing stepwise FWER controlling procedures,
including Bonferroni-Shaffer-based methods for testing logically related
hypotheses. However, this approach cannot be applied to dealing with the
directional errors problem.  The reason is that the approach can only exploit the logical relations among the tested hypotheses, whereas for developing powerful methods controlling directional errors, we need to exploit the special dependence relations of the test statistics as well as the logical relations of the tested hypotheses, as it is shown in the proof of Theorem 1.

\end{remark}

Although the upper bound of the $\text{FWER}_{\mathcal{F}_1}$ of
Procedure 1 is only slightly larger than $\alpha$,
this procedure cannot always control the
$\text{FWER}_{\mathcal{F}_1}$ at level $\alpha$ in the finite
samples. In the following, we present an example where the FWER of
the aforementioned procedure when testing $\mathcal{F}_1$ is above
$\alpha$ but of course below $\alpha/(1 - \alpha/n)$ as proved in
Theorem \ref{thm:1}.

\begin{example} \rm
Consider the special case of $\theta = (\theta_1, \ldots, \theta_n) \rightarrow 0$, thus $P_i \sim U(0,
1)$ for all $i \in I_0$. For Procedure 1, we have
\begin{eqnarray}\label{eq:10}
&& \text{FWER}_{\mathcal{F}_1}(\theta)  \nonumber \\
& = &  \sum_{r=0}^{n-1} {n
\choose r} \left (\frac{\alpha}{n} \right )^r \left [ \left ( 1 -
\frac{\alpha}{n} \right )^{n-r} - \left ( 1 - \frac{\alpha}{n} -
\frac{\alpha}{n -r} \right )^{n-r} \right ],
\end{eqnarray}
whose proof is given in the Appendix.
Through simple algebra calculation, we find out that
$\text{FWER}_{\mathcal{F}_1}(\theta) = \alpha + \frac{\alpha^2}{4} >
\alpha$ as $n=2$ and $\text{FWER}_{\mathcal{F}_1}(\theta) = \alpha +
\frac{\alpha^3}{108} > \alpha$ as $n=3$. Thus, the
Procedure 1 and thereby the usual Holm procedure with the
critical values $\alpha_i = \frac{\alpha}{n-i+1}, i = 1, \ldots, n$,
cannot always control the $\text{FWER}_{\mathcal{F}_1}$ at level
$\alpha$. \qed
\end{example}
It should be noted that in the above example,
 assumption A.2 is not used. This example shows that no matter whether or not
assumption A.2 holds, Procedure 1 cannot
control the FWER at level $\alpha$ in the finite samples.

\subsection{Holm-type stepdown procedure}

Consider a modified Holm procedure for testing $\mathcal{F}_1$ based on one-sided
$p$-values $P_i, i=1, \ldots, 2n$ defined in Section 2.3, which is described as follows.
\begin{procedure}
The stepdown procedure with the critical values $\alpha_i =
\frac{\alpha}{n - i + 1 + \alpha}, i = 1, \ldots, n$.
\end{procedure}
For any given parameter
vector $\theta = (\theta_1, \ldots, \theta_n)$, we have
\begin{eqnarray}\label{eq:14}
\text{FWER}_{\mathcal{F}_1}(\theta) \le \alpha,
\end{eqnarray}
whose proof is given in the appendix.
\begin{theorem}\label{thm:2}
Consider Procedure 3 defined as above.
Under assumptions A.1 - A.3, the procedure strongly controls the
$\text{FWER}_{\mathcal{F}_1}$ at level $\alpha$.
\end{theorem}

\begin{remark} \rm

It should be noted that if one directly uses the conventional Holm procedure
with the critical constants $\alpha_i = \alpha/(2n - i + 1), i = 1, \ldots, 2n$
for testing the $2n$ hypotheses, then the critical constants
corresponding to the first $n$ most significant hypotheses will be
always less than or equal to $\alpha/n$. However, for Procedure 3,
the critical constants corresponding to the first $n$
most significant hypotheses are generally much larger than
$\alpha/n$. The main reason why the Procedure 3 works well is that
the $2n$ tested hypotheses have some structural relationship: they
can be arranged as $n$ pairs of one true and one false null
hypotheses. For each pair of hypotheses, the sum of their
corresponding $p$-values is equal to one. Thus, for each pair of
hypotheses, when one hypothesis is significant, another one is
impossible to be significant. The newly introduced Procedure 3
has fully exploited the above facts and hence is more powerful than
the conventional Holm procedure.

\end{remark}

\begin{remark} \rm
It should be noted that when testing $n$ null hypotheses, the
critical constants of Procedure 3 are slightly less
than those of the usual Holm procedure, thus Procedure 3 can also
strongly control $\text{FWER}_{\mathcal{F}_2}$ at level $\alpha$.
Therefore, if we use separate analysis approach to test
$\mathcal{F}$ by applying separately Procedure 3 to test
$\mathcal{F}_1$ and $\mathcal{F}_2$ at level $\alpha/2$, then the
$\text{FWER}_{\mathcal{F}}$ is strongly controlled at level
$\alpha$.
\end{remark}

\section{Controlling the mdFWER under dependence}

In this section, we will discuss how to control the
$\text{FWER}_{\mathcal{F}_1}$ under three different types of
dependence: within- and between-block dependence, and positive
dependence.

\subsection{Controlling the $\text{FWER}_{\mathcal{F}_1}$ under block dependence}

Suppose that $\mathcal{F}_1 = \{H_1, \ldots, H_{2n}\}$ can be
organized as $b$ subfamilies $\mathcal{F}_{1i}, i = 1, \ldots, b$,
each of which have $n_i$ pairs of null hypotheses, $(H_j, H_{n+j})$,
with $\sum_{i=1}^b n_i = n$. Regarding the joint distribution of the
test statistics, except for positive dependence, the assumptions of two different types of block dependence are also invoked in the following sections. \vskip 5pt

\textbf{A.$3'$ } (\emph{Between-block dependence})  The test statistics corresponding
to the true null hypotheses within each subfamily $\mathcal{F}_{1i}, i = 1, \ldots, b$ are mutually
independent.  \vskip 5pt

\textbf{A.$3''$ } (\emph{Within-block dependence})  The test statistics corresponding
to the true null hypotheses between the subfamilies $\mathcal{F}_{1i}, i = 1, \ldots, b$ are mutually independent.  \vskip 5pt

By using Procedure 3, a method
for testing $\mathcal{F}_1$ can be constructed as follows:
\begin{procedure}
(Holm-type procedure under block dependence)
\begin{enumerate}
  \item For $i= 1, \ldots, b$, use Procedure 3 for
  testing $\mathcal{F}_{1i}$ at level $\beta_i = n_i \alpha/n$.
  \item Let $K_i$ be the corresponding set of rejected null
  hypotheses for testing $\mathcal{F}_{1i}$. Reject all null
  hypotheses in $\bigcup_{i=1}^b K_i$.
\end{enumerate}
\end{procedure}

Under the assumption of  between-block dependence, through Theorem
\ref{thm:2}, the FWER of Procedure 3 for testing
$\mathcal{F}_{1i}$, $\text{FWER}_{\mathcal{F}_{1i}}$, satisfies
$\text{FWER}_{\mathcal{F}_{1i}} \le n_i \alpha/n$. Thus, the overall
FWER of Procedure 4 for testing $\mathcal{F}_1$ satisfies
$$\text{FWER}_{\mathcal{F}_1} \le \sum_{i = 1}^b \text{FWER}_{\mathcal{F}_{1i}} \le \sum_{i = 1}^b \frac{n_i
\alpha}{n} = \alpha.$$
Therefore, we have the following result:

\begin{theorem}\label{thm:3}
Consider Procedure 4 defined as above. Under assumptions A.1, A.2 and A.$3'$, this procedure strongly
controls the
  $\text{FWER}_{\mathcal{F}_1}$ at level $\alpha$.
\end{theorem}

\begin{remark} \rm
When the number of subfamilies $b$ is equal to $n$, that is, each
subfamily has only one pair of hypotheses, Procedure 4 reduces
to a modified Bonferroni procedure with the critical constant
$\alpha/(n + \alpha)$, which strongly controls the
$\text{FWER}_{\mathcal{F}_1}$ under arbitrary dependence. When there
is only one subfamily, Procedure 4 reduces to Procedure 3, which strongly controls the
$\text{FWER}_{\mathcal{F}_1}$ under independence. Finally, we should
point out that the critical constants $\frac{n_i \alpha/n}{n_i - j +
1 + n_i \alpha/n}$ of the stepdown procedure used in
Procedure 4 are almost always larger than or equal to $\alpha/n$, which
implies that the method is generally more powerful than the usual
Bonferroni procedure with the critical constant $\alpha/n$.
\end{remark}

\begin{remark} \rm
When the test statistics corresponding to the above $b$ subfamilies
are within-block dependent rather than between-block dependent, we
can reorganize these $b$ subfamilies as $n_{\max}$ new subfamilies
such that the corresponding test statistics are between-block
dependent, where $n_{\max} = \max \{n_i: i = 1, \ldots, b \}.$ Then,
we can apply Procedure 4 to test $\mathcal{F}_1$
based on these reorganized subfamilies and it results in the
corresponding $\text{FWER}_{\mathcal{F}_1}$ is controlled at level
$\alpha$.
\end{remark}

\subsection{Controlling the $\text{FWER}_{\mathcal{F}}$ under positive dependence}

In this subsection, we discuss how to control the
$\text{FWER}_{\mathcal{F}}$ rather than
$\text{FWER}_{\mathcal{F}_1}$ under positive dependence. First,
reorganize $\mathcal{F}$ as two new subfamilies, $\mathcal{F}'_1 =
\{H_{i1}: i = 1, \ldots, n\}$ and $\mathcal{F}'_2 = \{H_{ij}: i = 1,
\ldots, n, j = 2, 3\}$. Thus, for each $i= 1, 2$, the test
statistics corresponding to the null hypotheses in $\mathcal{F}'_i$
are positively dependent (which is not the case for  $\mathcal{F}_1$, leading to the
current division into subfamilies).

Based on the conventional Hochberg procedure (Hochberg, 1988),
 which is the stepup procedure with  critical constants
$\alpha_i = \alpha/(n -i +1), i = 1, \ldots, n$ that strongly controls the
FWER at level $\alpha$ under positive dependence, a method for simultaneously
testing $\mathcal{F}$ can be constructed as follows:
\begin{procedure}
(Hochberg-type procedure under positive dependence)
\begin{enumerate}
  \item Use the Hochberg procedure to test $\mathcal{F}'_{1}$ at level $\alpha/2$.
  \item Use the Hochberg-type procedure with the critical constants \\ $\alpha_i = \frac{\alpha}{n-\lfloor (i+1)/2 \rfloor +1}, i = 1, \ldots 2n$, to test $\mathcal{F}'_{2}$ at level $\alpha/2$.
  \item For $i = 1, 2$, let $K_i$ be the corresponding set of rejected null
  hypotheses for testing $\mathcal{F}'_{i}$. Reject all null
  hypotheses in $K_1 \bigcup K_2$.
\end{enumerate}
\end{procedure}

Note that for $H_{i2}$ and $H_{i3}$, their corresponding $p$-values
are the same. Thus, when we apply the aforementioned Hochberg-type procedure in Procedure 5 to test $\mathcal{F}'_2$ at level $\alpha/2$, it is equivalent to apply the conventional Hochberg procedure with the critical constants $\alpha_i = \alpha/(n -i +1), i = 1, \ldots, n$ to test $H_{i2}$'s or
$H_{i3}$'s. Then, the corresponding
$\text{FWER}_{\mathcal{F}'_{2}}$ is controlled at level $\alpha/2$. (Of course, $\alpha$ could be split into $\beta$ and $\alpha - \beta$, but for simplicity $\beta = \alpha /2$.) Hence,
$$\text{FWER}_{\mathcal{F}} \le \text{FWER}_{\mathcal{F}'_{1}} + \text{FWER}_{\mathcal{F}'_{2}} \le
\alpha/2 + \alpha/2 = \alpha.$$

\begin{theorem}\label{thm:4}
Consider Procedure 5 defined as above. Under assumption A.1 and the
assumption of  positive dependence in the sense of (\ref{eq:3}), this procedure strongly controls
the
  $\text{FWER}_{\mathcal{F}}$ at level $\alpha$.
\end{theorem}

\section{Controlling the mixed directional FDR under independence and dependence}

In this section, we discuss how to control the
$\text{FDR}_{\mathcal{F}_1}$ under the same settings as in the last
two sections.

\subsection{On the $\text{FDR}_{\mathcal{F}_1}$ control under independence}

Consider the BH procedure (Benjamini and Hochberg, 1995) for testing $\mathcal{F}_1 = \{ H_1,
\ldots, H_{2n} \}$ based on one-sided $p$-values $P_i, i=1, \ldots, 2n$ defined in Section 2.3, which is described as follows.
\begin{procedure}
The stepup procedure with the critical values $\alpha_i = i\alpha/n, i = 1, \ldots, n$.
\end{procedure}
Note that $P_{n+i} = 1 - P_i$ for each $i = 1,
\ldots, n$;  thus, among the $2n$ corresponding $p$-values, there are
$n$ $p$-values larger than or equal to $0.5$. Therefore, for the
BH-type procedure, it is sufficient to only define its first $n$
critical constants while testing those $2n$ null hypotheses. Under
assumptions A.1 and A.3, for any given parameter vector $\theta =
(\theta_1, \ldots, \theta_n)$, we have
\begin{eqnarray}\label{eq:15}
\text{FDR}_{\mathcal{F}_1}(\theta)
& \le & \alpha,
\end{eqnarray}
whose proof is given in the Appendix. Therefore, the following conclusion holds.

\begin{theorem}\label{thm:5}
Consider Procedure 6 defined as above. Under assumptions A.1  and A.3,
the procedure strongly controls the $\text{FDR}_{\mathcal{F}_1}$ at
level $\alpha$.
\end{theorem}

\begin{remark}\rm
Note that assumption A.2 is not used.  In fact, the result holds without the parametric
model assumptions used in much of this paper. Indeed, all that is assumed is
the availability of  $p$-values $P_i$ for testing some parameter $\theta_i = 0$
and their independence.  Of course, we must have $P_{n+i} = 1- P_i$, but this
is a natural requirement when constructing two one-sided $p$-values.
\end{remark}

\begin{remark} \rm
When $\theta = 0$, the inequality in (\ref{eq:15}) becomes an equality. Thus
Procedure 6 cannot be improved in terms of its
critical values while maintaining the control of the
$\text{FDR}_{\mathcal{F}_1}$.
\end{remark}

\subsection{On the $\text{FDR}_{\mathcal{F}_1}$ control under between-block dependence}

Suppose that $\mathcal{F}_1 = \{H_1, \ldots, H_{2n}\}$ can be
organized as $b$ subfamilies $\mathcal{F}_{1i}, i = 1, \ldots, b$,
each of which have $n_i$ pairs of null hypotheses $(H_j, H_{n+j})$
with $\sum_{i=1}^b n_i = n$. Assume that the test statistics
corresponding to those subfamilies
satisfy the condition of between-block dependence.

By using Procedure 6, a method for
simultaneously testing $\mathcal{F}_1$ can be constructed as
follows:
\begin{procedure}
(BH-type procedure under between-block dependence)
\begin{enumerate}
  \item For each given $i= 1, \ldots, b$, use Procedure 6
  to test $\mathcal{F}_{1i}$ at level $n_i \alpha/n$.
  \item Let $K_i$ be the corresponding set of rejected null
  hypotheses for testing $\mathcal{F}_{1i}$. Reject all null
  hypotheses in $\bigcup_{i=1}^b K_i$.
\end{enumerate}
\end{procedure}

Under the assumption of  between-block dependence, through Theorem
5, the FDR of Procedure 6 for testing subfamily
$\mathcal{F}_{1i}$ at level $n_i \alpha/n$ satisfies $\text{FDR}_{\mathcal{F}_{1i}} \le n_i
\alpha/n$. Thus, the overall FDR of the above Procedure 7 for
testing $\mathcal{F}_1$ satisfies
$$\text{FDR}_{\mathcal{F}_1} \le \sum_{i = 1}^b \text{FDR}_{\mathcal{F}_{1i}} \le \sum_{i = 1}^b \frac{n_i
\alpha}{n} = \alpha.$$

\begin{theorem}\label{thm:6}
Consider Procedure 7 defined as above. Under assumptions A.1 and A.$3'$, this method strongly
controls the
  $\text{FDR}_{\mathcal{F}_1}$ at level $\alpha$.
\end{theorem}
Theorem \ref{thm:6} implies Theorem \ref{thm:5}. When $b=1$, it reduces to Theorem \ref{thm:5}.

\subsection{On the $\text{FDR}_{\mathcal{F}_1}$ control under within-block dependence}

Suppose that $\mathcal{F}_1 = \{H_1, \ldots, H_{2n}\}$ can be
organized as $b$ subfamilies $\mathcal{F}_{1i}, i = 1, \ldots, b$,
each of which have $n_i$ pairs of null hypotheses $(H_j, H_{n+j})$
with $\sum_{i=1}^b n_i = n$. Assume that the test statistics
corresponding to those subfamilies satisfy the condition of
within-block dependence. Note that there are exactly $n$ true null
hypotheses in $\mathcal{F}_1$, by exploiting the information in a
two-stage BH-type procedure introduced in Guo and Sarkar (2014), a
method for simultaneously testing $\mathcal{F}_1$ is constructed as
follows:
\begin{procedure}
(BH-type procedure under within-block dependence)
\begin{enumerate}
\item For $i = 1, \ldots, b$, let
$\widetilde{P}_i$ denote the smallest one among the $n_i$ pairs of
$p$-values corresponding to the $n_i$ pairs of null hypotheses in
$\mathcal{F}_{1i}$.
\item Order the smallest $p$-values $\widetilde{P}_i, i = 1, \ldots, b$
 as $\widetilde P_{(1)} \le \cdots \le \widetilde P_{(b)}$, and find $B = \max \{1 \le i \le b:
\tilde P_{(i)} \le i \alpha/n \}$.
\item In each subfamily $\mathcal{F}_{1i}$, reject those null hypotheses whose corresponding $p$-values
are less than or equal to $B\alpha/n$.
\end{enumerate}
\end{procedure}

By using the same arguments as in Guo and Sarkar (2012), we can show
that the above Procedure 8 strongly controls the
$\text{FDR}_{\mathcal{F}_1}$ at level $\alpha$. Therefore, we have
the following result.

\begin{theorem}\label{thm:7}
Consider Procedure 8 defined as above. Under assumptions A.1 and A.$3''$, this method strongly
controls the
  $\text{FDR}_{\mathcal{F}_1}$ at level $\alpha$.
\end{theorem}
Theorem \ref{thm:7} implies Theorem \ref{thm:5}. When $b=n$, it reduces to Theorem \ref{thm:5}.

\subsection{On the $\text{FDR}_{\mathcal{F}_1}$ control under positive dependence}

Suppose that the test statistics $T_i, i = 1, \ldots, n$ are
positively dependent in the sense of (\ref{eq:3}). Then, for each $j= 1, 2$,
the test statistics corresponding to the true null hypotheses $H_{ij}, i
= 1, \ldots, n$ are also positively dependent. For $j = 1, 2$, let
$\mathcal{F}_{1j} = \{H_{ij}, i = 1, \ldots, n\}$ and $n_{1j}$
denote the number of true nulls in $\mathcal{F}_{1j}$. Note that
there are exactly $n$ true null hypotheses in $\mathcal{F}_1 =
\mathcal{F}_{11}\bigcup\mathcal{F}_{12}$, thus $n_{11} + n_{12} =
n$. By using the similar idea due to Benjamini and Yekutieli (2005),
a method for testing $\mathcal{F}_1$ can be constructed as follows:
\begin{procedure}
(BH-type procedure under positive dependence)
\begin{enumerate}
  \item For $j= 1, 2$, use Procedure 6
  to test $\mathcal{F}_{1j}$ at level $\alpha$.
  \item Let $K_j$ be the corresponding set of rejected null
  hypotheses for testing $\mathcal{F}_{1j}$. Reject all null
  hypotheses in $K_1 \bigcup K_2$.
\end{enumerate}
\end{procedure}

By using the result in Benjamini and Yekutieli (2001) and Sarkar
(2002) on the FDR control of the BH procedure under positive
dependence, we have
$$\text{FDR}_{\mathcal{F}_1} \le \text{FDR}_{\mathcal{F}_{11}} + \text{FDR}_{\mathcal{F}_{12}} \le
\frac{n_{11} \alpha}{n} + \frac{n_{12} \alpha}{n} = \alpha.$$ The
equality follows from the fact that $n_{11} + n_{12} = n$.

\begin{theorem}\label{thm:8}
Consider Procedure 9 defined as above. Under assumption A.1 and the assumption of  positive dependence in the sense of (\ref{eq:3}), this method
strongly controls the
  $\text{FDR}_{\mathcal{F}_1}$ at level $\alpha$.
\end{theorem}

\section{Concluding remarks}

In this paper, several approaches, methods, and results are presented addressing the
multiple testing problem of accounting for both type 1 and type 3 errors.
Many of the results required the assumption of independence, which is quite strong,
though we have weakened this assumption as well.   The problem
of directional error control has proven to be quite challenging, and though we
do not consider the dependent case more fully, it is hoped to consider this important
problem in future work.

\vskip 10pt

\section*{Acknowledgements}
The research of the first author was supported in part
by NSF Grant DMS-1006021 and DMS-1309162 and the
research of the second author was supported in part
by NSF Grant DMS-0707085.

\vskip 10pt

\section*{Appendix: Proofs} \vskip 5pt

\subsection*{A.1. Proof of Lemma 1}

Since the family of densities $f_{i, \theta}(\cdot)$ satisfies the assumption of MLR,
we have that,  for any given $\theta_1 > \theta_0$ and $x_1 > x_0$,
\begin{equation}\label{eq:old19}
\frac{f_{i, \theta_1}(x_1)}{f_{i, \theta_0}(x_1)} \ge
\frac{f_{i, \theta_1}(x_0)}{f_{i, \theta_0}(x_0)}.
\end{equation}

By multiplying both sides of (\ref{eq:old19}) by $f_{i, \theta_0} ( x_0 )$ and then integrating
over $x_0$ from $- \infty $ to $x_1$, one obtains
\begin{equation}\label{eq:old20}
\frac{f_{i, \theta_1}(x)}{f_{i, \theta_0}(x)} \ge
\frac{F_{i, \theta_1}(x)}{F_{i, \theta_0}(x)}~.
\end{equation}
Similarly, one obtains
\begin{equation}\label{eq:new20}
 \frac{1 -
F_{i, \theta_1}(x)}{1 - F_{i, \theta_0}(x)}  \ge
\frac{f_{i, \theta_1}(x)}{f_{i, \theta_0}(x)}.
\end{equation}

Consider the functions $G_1(x) =
\frac{F_{i, \theta_1}(x)}{F_{i, \theta_0}(x)}$ and $G_2(x) = \frac{1 -
F_{i, \theta_1}(x)}{1 - F_{i, \theta_0}(x)}$. It is easy to check by using
(\ref{eq:old20})  and (\ref{eq:new20})  that $G'_1(x) \ge 0$ and $G'_2(x) \ge 0$. Then, $G_1(x)$ and
$G_2(x)$ are both non-decreasing in $x$.
First, assume $\theta_i > 0$, so  that $P_i =  F_{i, 0} ( T_i )$.  Thus, for any $0 \le p < p'
\le 1$,
\begin{eqnarray}\label{eq:old21}
\text{Pr}_{\theta_i} \left \{ P_i \le p \big | P_i \le p' \right \}
= \frac{F_{i, \theta_i}(x)}{F_{i, \theta_i}(x')} \le
\frac{F_{i, 0}(x)}{F_{i, 0}(x')} = \frac{p}{p'} = \text{Pr}_{\theta_i = 0}
\left \{ P_i \le p \big | P_i \le p' \right \},
\end{eqnarray}
where $x = F_0^{-1}(p)$ and $x' = F_0^{-1}(p')$. In (\ref{eq:old21}), the
inequality follows from the fact that $G_1(x)$ is non-decreasing in
$x$ and the second equality follows from assumption A.1. By using
 similar arguments, we can prove that (\ref{eq:old21}) also holds when
$\theta_i < 0$. Hence, the desired result follows. ~\qed

\vskip 10pt

\subsection*{A.2. Proof of (\ref{eq:6})}

Throughout the Appendix, the following
notation will be used.  Given any index set of false null hypotheses, $S
\subset I_1$,  define $\overline{S} = I_1 \backslash S$,
$S_{-n} = \{i \in I_0: n+i \in S \}$, and
$\overline{S}_{-n} = \{i \in I_0: n+i \in \overline{S} \}$.
It is easy to see that $|S| = |S_{-n}|$ and $|\overline{S}| =
|\overline{S}_{-n}|.$

Consider Procedure 1 for testing $\mathcal{F}_1$. Let $R_{11}$ be the index set of rejected false null hypotheses
at the first stage, $R_{10}$ be the index set of true null hypotheses
for which the corresponding $p$-values less than $1 - \alpha/n$, and
$R_{10}^{(-j)}$ be the index set of true null hypotheses excluding
$H_j$ for which the corresponding $p$-values less than $1 -
\alpha/n$, that is, $R_{11} = \{i \in I_1: P_i \le \alpha/n \},
R_{10} = \{i \in I_0: P_i < 1 - \alpha/n \}$, and $R_{10}^{(-j)} =
\{i \in I_0 \backslash \{ j \}: P_i < 1 - \alpha/n \} = R_{10}
\backslash \{ j \}$.

Let $\widehat P_{(1)}^{I_0} $ be the minimum $p$-value corresponding to
the true null hypotheses with indices in $I_0$, for any given parameter
vector $\theta = (\theta_1, \ldots, \theta_n)$, we have
\begin{eqnarray}\label{eq:4_2}
& & \text{FWER}_{\mathcal{F}_1}(\theta)  =  \sum_{S \subset I_1} \text{Pr}_{\theta}
\left \{R_{11} = S, \widehat P_{(1)}^{I_0} \le \frac{\alpha}{n - |S|}
\right \} \nonumber \\
& = & \sum_{S \subset I_1} \text{Pr}_{\theta} \left \{P_i \le
\frac{\alpha}{n} \text{ for~all } i \in S, P_i > \frac{\alpha}{n}
\text{ for~all } i \in \overline{S}, \widehat P_{(1)}^{I_0} \le
\frac{\alpha}{n - |S|}
\right \} \nonumber \\
& = & \sum_{S_{-n} \subset I_0} \text{Pr}_{\theta} \left \{P_i \ge 1-
\frac{\alpha}{n} \text{ for~all } i \in S_{-n}, \right . \nonumber \\
& & \hspace{10em} \left . P_i <
1-\frac{\alpha}{n} \text{ for~all } i \in \overline{S}_{-n}, \widehat
P_{(1)}^{\overline{S}_{-n}} \le \frac{\alpha}{n - |S_{-n}|}
\right \} \nonumber \\
& \le & \sum_{S_{-n} \subset I_0} \sum_{j \in \overline{S}_{-n}}
\text{Pr}_{\theta} \left \{P_i \ge 1- \frac{\alpha}{n} \text{ for~all } i \in
S_{-n}, \right . \nonumber \\
& & \hspace{10em} \left . P_i < 1-\frac{\alpha}{n} \text{ for~all } i \in \overline{S}_{-n},
P_j \le \frac{\alpha}{n - |S_{-n}|}
\right \} \nonumber \\
& = & \sum_{S_{-n} \subset I_0} \sum_{j \in \overline{S}_{-n}}
\text{Pr}_{\theta} \left \{{R}_{10}^{(-j)} = \overline{S}_{-n}
\backslash \{ j \}, P_j \le \frac{\alpha}{n - |S_{-n}|} \right \}.
\end{eqnarray}
The inequality follows from the Bonferroni inequality.

Note that for $j \in \overline{S}_{-n}$,
\begin{eqnarray}\label{eq:5_2}
& & \text{Pr}_{\theta} \left \{{R}_{10}^{(-j)} = \overline{S}_{-n}
\backslash \{ j \}, P_j \le \frac{\alpha}{n - |S_{-n}|}
\right \} \nonumber \\
& = & \text{Pr}_{\theta} \left \{{R}_{10}^{(-j)} = \overline{S}_{-n}
\backslash \{ j \}, P_j \le \frac{\alpha}{n - |S_{-n}|}, P_j \le 1 -
\frac{\alpha}{n}
\right \} \nonumber \\
& = & \text{Pr}_{\theta} \left \{{R}_{10}^{(-j)} = \overline{S}_{-n}
\backslash \{ j \}, P_j \le \frac{\alpha}{n - |S_{-n}|} \Big | P_j
\le 1 - \frac{\alpha}{n} \right \} \nonumber \\
& & \hspace{17em} \times ~ \text{Pr}_{\theta_j} \left \{ P_j
\le 1 - \frac{\alpha}{n}
\right \} \nonumber \\
& = & \text{Pr}_{\theta^{(-j)}} \left \{{R}_{10}^{(-j)} =
\overline{S}_{-n} \backslash \{ j \} \right \} \text{Pr}_{\theta_j}
\left \{ P_j \le 1 - \frac{\alpha}{n} \right \} \nonumber \\
& & \hspace{11em} \times ~ \text{Pr}_{\theta_j} \left \{ P_j \le \frac{\alpha}{n - |S_{-n}|} \Big | P_j \le 1 - \frac{\alpha}{n} \right \} \nonumber \\
& \le & \text{Pr}_{\theta} \left \{{R}_{10} = \overline{S}_{-n} \right \} \text{Pr}_{\theta_j = 0} \left \{ P_j \le \frac{\alpha}{n - |S_{-n}|} \Big | P_j \le 1 - \frac{\alpha}{n} \right \} \nonumber \\
& = & \text{Pr}_{\theta} \left \{{R}_{10} = \overline{S}_{-n} \right
\} \frac{\alpha}{n - |S_{-n}|} \frac{1}{1 - \alpha/n} ,
\end{eqnarray}
where $\theta^{(-j)} = (\theta_1, \ldots, \theta_{j-1},
\theta_{j+1}, \ldots, \theta_n)$. Here, the third equality follows
from assumption A.3 and the fourth follows from assumption A.1 under
which $P_j \sim U(0, 1)$ when $\theta_j = 0$. For the inequality,
the first term of its right-hand side follows from assumption A.3
under which the first two terms of the left-hand side match up, and
the second one of its right-hand side follows from assumption A.2.

Applying (\ref{eq:5_2}) to (\ref{eq:4_2}), we have
\begin{eqnarray}\label{eq:6_2}
\text{FWER}_{\mathcal{F}_1}(\theta) & \le & \sum_{S_{-n} \subset
I_0} \sum_{j \in \overline{S}_{-n}} \text{Pr}_{\theta} \left
\{{R}_{10} = \overline{S}_{-n} \right \} \frac{\alpha}{n - |S_{-n}|}
\frac{1}{1
- \alpha/n} \nonumber \\
& = & \sum_{S_{-n} \subset I_0} \text{Pr}_{\theta} \left \{ {R}_{10} = \overline{S}_{-n} \right \} \frac{\alpha}{1 - \alpha/n} \nonumber \\
& \le & \frac{\alpha}{1 - \alpha/n}.
\end{eqnarray}
Hence, the desired result follows. \qed

\vskip 10pt

\subsection*{A.3. Proof of (\ref{eq:10})}
By using the third equality of (\ref{eq:4_2}), we have
\begin{eqnarray}\label{eq:8}
& & \text{FWER}_{\mathcal{F}_1}(\theta) \nonumber \\
& = & \sum_{r=0}^{n-1} \sum_{\substack{S_{-n} \subset I_0 \\
|S_{-n}| = r}} \text{Pr}_{\theta} \left \{P_i \ge 1-
\frac{\alpha}{n} \text{ for } i \in S_{-n}, \right . \nonumber \\
& & \hspace{13em} \left .
P_i <
1-\frac{\alpha}{n} \text{ for } i \in \overline{S}_{-n},
 \widehat P_{(1)}^{\overline{S}_{-n}} \le \frac{\alpha}{n - r}
\right \}  \nonumber \\
& = & \sum_{r=0}^{n-1} \sum_{\substack{S_{-n} \subset I_0 \\
|S_{-n}| = r}} \text{Pr}_{\theta} \left \{P_i \ge 1-
\frac{\alpha}{n} \text{ for } i \in S_{-n} \right \} \nonumber \\
& & \hspace{9em} \times ~ \text{Pr}_{\theta} \left \{ P_i < 1-\frac{\alpha}{n} \text{ for } i
\in \overline{S}_{-n}, \widehat P_{(1)}^{\overline{S}_{-n}} \le
\frac{\alpha}{n - r} \right \} \nonumber \\
& = & \sum_{r=0}^{n-1} \sum_{\substack{S_{-n} \subset I_0 \\
|S_{-n}| = r}} \left (\frac{\alpha}{n} \right )^r \text{Pr}_{\theta}
\left \{ P_i < 1-\frac{\alpha}{n} \text{ for } i \in
\overline{S}_{-n}, \widehat P_{(1)}^{\overline{S}_{-n}} \le
\frac{\alpha}{n - r} \right \}.
\end{eqnarray}
In the above special case with $|S_{-n}| = r$, we have
\begin{eqnarray}\label{eq:9}
& & \text{Pr}_{\theta} \left \{P_i < 1-\frac{\alpha}{n} \text{
for~all } i \in \overline{S}_{-n}, \widehat
P_{(1)}^{\overline{S}_{-n}} \le \frac{\alpha}{n - r}
\right \} \nonumber \\
& = & \text{Pr}_{\theta} \left \{\widehat
P_{(1)}^{\overline{S}_{-n}} \le \frac{\alpha}{n - r} \big |
P_i < 1-\frac{\alpha}{n} \text{ for~all } i \in
\overline{S}_{-n}  \right \} \nonumber \\
& & \hspace{11em} \times ~ \text{Pr}_{\theta} \left \{ P_i
< 1-\frac{\alpha}{n} \text{ for~all
} i \in \overline{S}_{-n}  \right \} \nonumber \\
& = & \left [ 1 -  \text{Pr}_{\theta} \left \{\widehat
P_{(1)}^{\overline{S}_{-n}} > \frac{\alpha}{n - r} \big | P_i <
1-\frac{\alpha}{n} \text{ for~all } i \in \overline{S}_{-n}  \right
\} \right ] \nonumber \\
& & \hspace{15em} \times ~ \prod_{ i \in \overline{S}_{-n}} \text{Pr}_{\theta}
\left \{ P_i < 1-\frac{\alpha}{n} \right \}  \nonumber \\
& = & \left [ 1 -  \prod_{ i \in \overline{S}_{-n}}
\text{Pr}_{\theta} \left \{P_i > \frac{\alpha}{n - r} \big | P_i <
1-\frac{\alpha}{n} \right \} \right ] \left ( 1 - \frac{\alpha}{n}
\right )^{n-r} \nonumber \\
& = & \left [ 1 -  \left (1 - \frac{\frac{\alpha}{n-r}}{1 -
\frac{\alpha}{n}}
\right )^{n-r} \right ] \left ( 1 - \frac{\alpha}{n} \right )^{n-r} \nonumber \\
& = & \left ( 1 - \frac{\alpha}{n} \right )^{n-r} - \left ( 1 -
\frac{\alpha}{n} - \frac{\alpha}{n -r} \right )^{n-r}.
\end{eqnarray}
Here, the second and third equalities follow from assumption A.3.
Apply (\ref{eq:9}) into (\ref{eq:8}), we have
\begin{eqnarray}
&& \text{FWER}_{\mathcal{F}_1}(\theta)  \nonumber \\
& = &  \sum_{r=0}^{n-1} {n
\choose r} \left (\frac{\alpha}{n} \right )^r \left [ \left ( 1 -
\frac{\alpha}{n} \right )^{n-r} - \left ( 1 - \frac{\alpha}{n} -
\frac{\alpha}{n -r} \right )^{n-r} \right ], \nonumber
\end{eqnarray}
the desired result. \qed

\vskip 10pt

\subsection*{A.4. Proof of (\ref{eq:14})}

Consider Procedure 3 for testing $\mathcal{F}_1$. Let $\widehat q_{(1)} \le \ldots \le \widehat
q_{(n)}$ denote the ordered false null $p$-values. Define $J = \max
\{j: \widehat q_{(i)} \le \alpha_i, ~\forall i \le j \}$, provided
this maximum exists; otherwise, let $J = 0$. Let $K$ denote the
index set of the $J$ rejected false null hypotheses when applying
the stepdown procedure to simultaneously test the $n$ false null
hypotheses $H_{n+1}, \ldots, H_{2n}$, and $E_1$ denote the event of
at least one falsely rejected hypothesis when applying the same
procedure to simultaneously test $H_1, \ldots, H_{2n}$. It should be
noted that if $J = n$, then no true null hypotheses are falsely
rejected when testing $\mathcal{F}_1$.  Thus,
\begin{eqnarray}\label{eq:11}
E_1 & = & \bigcup_{j =0}^{n-1} \left \{ J = j, \widehat
P_{(1)}^{I_0}
\le \alpha_{j+1} \right \} \nonumber \\
& = & \bigcup_{j=0}^{n-1} \bigcup_{\substack{ S \subset I_1 \\
|S| = j}} \left \{ K = S, \widehat P_{(1)}^{I_0}
\le \alpha_{j+1} \right \} \nonumber \\
& = & \bigcup_{S \subset I_1} \left \{ K = S, \widehat
P_{(1)}^{\overline{S}_{-n}} \le \alpha_{|S|+1} \right \}.
\end{eqnarray}
For any given parameter vector $\theta =
(\theta_1, \ldots, \theta_n)$, we have
\begin{eqnarray}\label{eq:12}
\text{FWER}_{\mathcal{F}_1}(\theta) & = & \text{Pr}_{\theta}(E_1)
\nonumber \\
& = & \sum_{S \subset I_1} \text{Pr}_{\theta} \left \{ K = S,
\widehat P_{(1)}^{\overline{S}_{-n}} \le
\alpha_{|S|+1} \right \} \nonumber \\
& \le & \sum_{S \subset I_1} \sum_{j \in \overline{S}_{-n}}
\text{Pr}_{\theta} \left \{ K = S,
P_j \le \alpha_{|S|+1} \right \}  \nonumber \\
& = & \sum_{S \subset I_1} \sum_{j \in \overline{S}}
\text{Pr}_{\theta} \left \{ K^{\{ -j \}} = S, P_j \ge 1 -
\alpha_{|S|+1} \right \},
\end{eqnarray}
where $K^{\{ -j \}}$ is the index set of rejected false null
hypotheses by using the stepdown procedure with the critical
constants $\alpha_i = \frac{\alpha}{n - i + 1 + \alpha}, i = 1,
\ldots, n-1$ to simultaneously test the $n-1$ false null hypotheses
$H_{n+1}, \ldots, H_{2n}$ excluding $H_{j}$ with $j \in I_1$.

By using the similar argument lines as in (\ref{eq:5_2}), we have
\begin{eqnarray}\label{eq:13}
& & \text{Pr}_{\theta} \left \{ K^{\{ -j \}} = S, P_j
\ge 1 - \alpha_{|S|+1} \right \} \nonumber \\
& = & \text{Pr}_{\theta} \left \{ K^{\{ -j \}} = S, P_j
> \alpha_{|S|+1} \right \} \text{Pr}_{\theta_j} \left \{ P_j
\ge 1 - \alpha_{|S|+1} \Big | P_j
> \alpha_{|S|+1} \right \} \nonumber \\
& \le & \text{Pr}_{\theta} \left \{K = S \right \}
\text{Pr}_{\theta_j = 0} \left \{ P_j \ge 1 - \alpha_{|S|+1} \Big |
P_j
> \alpha_{|S|+1} \right \} \nonumber \\
& = & \text{Pr}_{\theta} \left \{K = S \right \}
\alpha_{|S|+1} \frac{1}{1 - \alpha_{|S|+1}} \nonumber \\
& = & \text{Pr}_{\theta} \left \{K = S \right \} \frac{\alpha}{n -
|S|}.
\end{eqnarray}
Applying (\ref{eq:13}) to (\ref{eq:12}), we have
\begin{eqnarray}
\text{FWER}_{\mathcal{F}_1}(\theta) & \le & \alpha
\sum_{S \subset I_1} \text{Pr}_{\theta} \left \{ K = S \right \}  \le \alpha, \nonumber
\end{eqnarray}
the desired result.  ~\qed

\vskip 10pt

\subsection*{A.5. Proof of (\ref{eq:15})}
Consider Procedure 6 for testing $\mathcal{F}_1$. Note that under assumptions A.1 and A.3, for any given parameter vector $\theta =
(\theta_1, \ldots, \theta_n)$, we have
\begin{eqnarray}
\text{FDR}_{\mathcal{F}_1}(\theta) &=& E \left \{ \frac{V_1}{R_1 \vee 1} \right \} \nonumber \\
& = & \sum_{i=1}^n E_{\theta} \left \{ \frac{I \{ H_i~{\rm rejected} \}  }{R_1 \vee 1} \right \} \nonumber \\
& = & \sum_{i=1}^n \sum_{r=1}^n \frac{1}{r} P_{\theta} \{  R_1 = r,~H_i~{\rm rejected} \} \nonumber \\
& = & \sum_{i=1}^n \sum_{r=1}^n \frac{1}{r} \text{Pr}_{\theta} \left
(R_1 = r,
P_i \le \frac{r}{n} \alpha \right ) \nonumber \\
& = & \sum_{i=1}^n \sum_{r=1}^n \frac{1}{r} \text{Pr}_{\theta} \left
(R_1^{\{-i, -(n+i)\}} = r-1,
P_i \le \frac{r}{n} \alpha \right ) \nonumber \\
& \le & \sum_{i=1}^n \sum_{r=1}^n \frac{\alpha}{n}
\text{Pr}_{\theta} \left
(R_1^{\{-i, -(n+i)\}} = r-1 \right ) \nonumber \\
& = & \alpha. \nonumber
\end{eqnarray}
Here, $R_1^{\{-i, -(n+i)\}}$ is the number of rejected null
hypotheses by using the stepup procedure with the critical values
$j\alpha/n, j = 2, \ldots, n$ to simultaneously test the
$2(n-1)$ null hypotheses $H_1, \ldots, H_{2n}$ excluding the pair of
null hypotheses $(H_{i}, H_{n+i})$. The inequality follows
from assumptions A.1 and A.3 and the fact that $P_{n+i} = 1 - P_i$. ~\qed

\bigskip

\bibliographystyle{elsarticle-harv}

\end{document}